\documentclass[12pt]{amsart}
\usepackage{amssymb,amsmath,amsfonts,latexsym,setspace}
\usepackage{bm}
\newcommand{\newsection}[1]{\setcounter{equation}{0} \section{#1}}
\newcommand{\bea}{\begin{eqnarray}}
\newcommand{\eea}{\end{eqnarray}}

\newcommand{\cle}{\mathcal{E}}

\newcommand{\clh}{\mathcal{H}}
\newcommand{\clk}{\mathcal{K}}
\newcommand{\cll}{\mathcal{L}}

\newcommand{\clw}{\mathcal{W}}

\newcommand{\raro}{\rightarrow}

\def \qed {\hfill \vrule height6pt width 6pt depth 0pt}
\def\textmatrix#1&#2\\#3&#4\\{\bigl({#1 \atop #3}\ {#2 \atop #4}\bigr)}
\def\dispmatrix#1&#2\\#3&#4\\{\left({#1 \atop #3}\ {#2 \atop #4}\right)}
\newcommand{\be}{\begin{equation}}
\newcommand{\ee}{\end{equation}}
\newcommand{\ben}{\begin{eqnarray*}}
\newcommand{\een}{\end{eqnarray*}}

\newcommand{\NI}{\noindent}

\newcommand{\bi}{\begin{itemize}}
\newcommand{\ei}{\end{itemize}}

\newtheorem{Theorem}{\sc Theorem}[section]
\newtheorem{Lemma}[Theorem]{\sc Lemma}
\newtheorem{Proposition}[Theorem]{\sc Proposition}
\newtheorem{Corollary}[Theorem]{\sc Corollary}
\newtheorem{Definition}[Theorem]{\sc Definition}
\newtheorem{Example}[Theorem]{\sc Example}
\newtheorem{Remark}[Theorem]{\sc Remark}
\newtheorem{Note}[Theorem]{\sc Note}
\newtheorem{Question}{\sc Question}
\newtheorem{ass}[Theorem]{\sc Assumption}
\newcommand{\bt}{\begin{Theorem}}
\def\beginlem{\begin{Lemma}}
\def\beginprop{\begin{Proposition}}
\def\begincor{\begin{Corollary}}
\def\begindef{\begin{Definition}}
\def\beginexamp{\begin{Example}}
\def\beginrem{\begin{Remark}}
\def\beginq{\begin{Question}}
\def\beginass{\begin{ass}}
\def\beginnote{\begin{Note}}
\newcommand{\et}{\end{Theorem}}
\def\endlem{\end{Lemma}}
\def\endprop{\end{Proposition}}
\def\endcor{\end{Corollary}}
\def\enddef{\end{Definition}}
\def\endexamp{\end{Example}}
\def\endrem{\end{Remark}}
\def\endq{\end{Question}}
\def\endass{\end{ass}}
\def\endnote{\end{Note}}

\begin{document}

\title[Wold decomposition for doubly commuting isometries]{Wold decomposition for doubly commuting isometries}

\author[Jaydeb Sarkar]{Jaydeb Sarkar}
\address{Indian Statistical Institute, Statistics and Mathematics Unit, 8th Mile, Mysore Road, Bangalore, 560059, India}
\email{jay@isibang.ac.in, jaydeb@gmail.com}

\subjclass[2010]{47A13, 47A15, 47A20, 47L99 }


\keywords{Commuting isometries, Wold decomposition, Doubly
commutativity}

\begin{abstract} In this paper, we obtain a complete description of
the class of $n$-tuples ($n \geq 2$) of doubly commuting isometries.
In particular, we present a several variables analogue of the Wold
decomposition for isometries on Hilbert spaces. Our main result is a
generalization of M. Slocinski's Wold-type decomposition of a pair
of doubly commuting isometries.
\end{abstract}

\maketitle

\newsection{Introduction}

Let $V$ be an isometry on a Hilbert space $\clh$, that is, $V^* V =
I_{\clh}$. A closed subspace $\clw \subseteq \clh$ is said to be
\textit{wandering subspace} for $V$ if $V^k \clw \perp V^l \clw$ for
all $k, l \in \mathbb{N}$ with $k \neq l$, or equivalently, if $V^m
\clw \perp \clw$ for all $m \geq 1$. An isometry $V$ on $\clh$ is
said to be a \textit{unilateral shift} or \textit{shift} if \[\clh =
\mathop{\bigoplus}_{m \geq 0} V^m \clw,\]for some wandering subspace
$\clw$ for $V$.

For a shift $V$ on $\clh$ with a wandering subspace $\clw$ we have
\[\clh \ominus V \clh = \mathop{\bigoplus}_{m \geq 0} V^m
\clw \ominus V(\mathop{\bigoplus}_{m \geq 0} V^m \clw) =
\mathop{\bigoplus}_{m \geq 0} V^m \clw \ominus \mathop{\bigoplus}_{m
\geq 1} V^m \clw = \clw.\]In other words, the wandering subspace of
a shift is unique and is given by $\clw = \clh \ominus V \clh$. The
dimension of the wandering subspace of a shift is called the
\textit{multiplicity} of the shift.

The classical Wold decomposition theorem (\cite{W}, see also page 3
in \cite{NF}) states that every isometry on a Hilbert space is
either a shift, or a unitary, or a direct sum of shift and unitary
(see Theorem \ref{Wold}). The Wold decomposition theorem, or results
analogous of the Wold decomposition theorem, plays an important role
in many areas of operator algebras and operator theory, including
the invariant subspace problem for Hilbert spaces of holomorphic
functions (cf. \cite{KMN}, \cite{KM}, \cite{M}, \cite{SZ},
\cite{Sh}).

The natural question then becomes: Let $ n \geq 2$ and $V = (V_1,
\ldots, V_n)$ be an $n$-tuple of commuting isometries on a Hilbert
space $\clh$. Does there exists a Wold-type decomposition for $V$?
What happens if we have a family of isometries?

Several interesting results have been obtained in many directions.
For instance, in \cite{Su} Suciu developed a structure theory for
semigroup of isometries (see also \cite{GG}, \cite{HL}, \cite{KO}).
However, to get a more precise result it will presumably always be
necessary to make additional assumptions on the family of operators.
Before proceeding, we recall the definition of the doubly commuting
isometries.

Let $V = (V_1, \ldots, V_n)$ be an $n$-tuple $(n \geq 2$) of
commuting isometries on $\clh$. Then $V$ is said to \textit{doubly
commute} if\[V_i V^*_j = V^*_j V_i,\]for all $1 \leq i < j \leq n$.

The simplest example of an $n$-tuple of doubly commuting isometries
is the tuple of multiplication operators $(M_{z_1}, \ldots,
M_{z_n})$ by the coordinate functions on the Hardy space
$H^2(\mathbb{D}^n)$ over the polydisc $\mathbb{D}^n$ ($n \geq 2$).

In \cite{S}, M. Slocinski obtained an analogous result of Wold
decomposition theorem for a pair of doubly commuting isometries.

\begin{Theorem} \label{Sl} \textsf{(M. Slocinski)} Let $V =
(V_1, V_2)$ be a pair of doubly commuting isometries on a Hilbert
space $\clh$. Then there exists a unique decomposition
\[\clh = \clh_{ss} \oplus \clh_{su} \oplus \clh_{us} \oplus
\clh_{uu},\]where $\clh_{ij}$ are joint $V$-reducing subspace of
$\clh$ for all $i, j = s, u$. Moreover, $V_1$ on ${\clh_{i, j}}$ is
a shift if $i = s$ and unitary if $i = u$ and that $V_2$ is a shift
if $j =s$ and unitary if $j = u$.
\end{Theorem}

We refer to \cite{KO} for a new proof of Slocinski's result (see
also \cite{CPZ}, \cite{KMN}, \cite{KM}).

Slocinski's Wold-type decomposition does not hold for general tuples
of commuting isometries (cf. Example 1 in \cite{S}). However, if $V$
is an $n$-tuple of commuting isometries such that \[\dim \ker
\Big(\mathop{\prod}_{i=1}^n V_i^*\Big) < \infty,\] then $V$ admits a
Wold-type decomposition (see Theorem 2.4 in \cite{BKS}).

In this paper, we obtain a Wold-type decomposition for tuples of
doubly commuting isometries. We extend the ideas of M. Slocinski on
the Wold-type decomposition for a pair of isometries to the
multivariable case ($n \geq 2$). Our approach is simple and based on
the classical Wold decomposition for a single isometry. Moreover,
our method yields a new proof of Slocinski's result for the base
case $n = 2$. In addition, we obtain an explicit description of the
closed subspaces in the orthogonal decomposition of the Hilbert
space (see the equality (\ref{W2})).

The paper is organized as follows. In Section 2, we set up notations
and definitions and establish some preliminary results. In Section
3, we prove our main result and some of its consequences.

\NI\textsf{Note added in proof:} After this work was completed, we
became aware that generalization of Slocinski's results to
$n$-tuples of doubly commuting isometries has been obtained
independently by Timotin \cite{T}, Gaspar and Suciu \cite{GS} and in
the context of $C^*$-correspondence, by Skalski and Zacharias
\cite{SZ}. However, our main theorem, Theorem \ref{dc-w}, is
stronger in the sense that the orthogonal decomposition in
(\ref{W1}) works for any $m \in \{2, \ldots, n\}$ (with $2 < n$).
Also our approach yields an explicit representation of the closed
subspaces in the orthogonal decomposition (see the equality
(\ref{W2})) and a new and simpler proof of earlier generalizations
of Slocinski's result. It is a question of general interest whether,
the present method, with explicit description of the closed
subspaces (\ref{W2}) in the Wold decomposition (\ref{W1}) yields
more structure theory results in the study of $C^*$-correspondences.

\newsection{Preparatory results}

In this section we recall the Wold decomposition theorem and present
some elementary facts concerning doubly commuting isometries.

We begin with the Wold decomposition of an isometry.

\begin{Theorem}\label{Wold}\textsf{(H. Wold)}
Let $V$ be an isometry on $\clh$. Then $\clh$ admits a unique
decomposition $\clh = \clh_s \oplus \clh_u$, where $\clh_s$ and
$\clh_u$ are $V$-reducing subspaces of $\clh$ and $V|_{\clh_s}$ is a
shift and $V|_{\clh_u}$ is unitary. Moreover,
\[\clh_s = \mathop{\bigoplus}_{m=0}^\infty V^m \clw \quad \quad
\mbox{and} \quad \quad \clh_u = \bigcap_{m=0}^\infty V^m
\clh,\]where $\clw = \mbox{ran}(I - V V^*) = \mbox{ker} V^*$ is the
wandering subspace for $V$.
\end{Theorem}
\NI\textsf{Proof.} Let $\clw = \mbox{ran}(I - V V^*)$ be the
wandering subspace for $V$ and \[\clh_s : =
\mathop{\bigoplus}_{m=0}^\infty V^m \clw.\] Consequently, $\clh_s$
is a $V$-reducing subspace of $\clh$ and that $V|_{\clh_s}$ is an
isometry. Furthermore,\[\clh_u := \clh_s^\perp =
\big(\mathop{\bigoplus}_{m=0}^\infty V^m \clw \big)^\perp =
\mathop{\bigcap}_{m=0}^\infty (V^m \clw)^\perp.\]We observe now that
$I - V V^*$ is an orthogonal projection, hence $V^l (I - V V^*)
V^{*l}$ is also an orthogonal projection and \[V^l (I - V V^*)
V^{*l} = (V^l (I - V V^*)) (V^l (I - V V^*))^*,\]for all $l \geq 0$.
Thus we obtain \[\mbox{ran} V^l (I - V V^*) = \mbox{ran} \big((V^l
(I - V V^*)) (V^l (I - V V^*))^*\big) = \mbox{ran} V^l (I - V V^*)
V^{*l},\]and hence \[\begin{split} (V^l \clw)^\perp & = (V^l
\mbox{ran}(I - V V^*))^\perp = (\mbox{ran} V^l (I - V V^*))^\perp \\
& = (\mbox{ran} V^l (I - V V^*) V^{*l})^\perp = \mbox{ran} (I - V^l
(I - V V^*) V^{*l}) \\ & = \mbox{ran}[ (I - V^l V^{*l}) \oplus
V^{l+1} V^{* \,l+1}] = \mbox{ran}(I - V^l V^{*l}) \oplus \mbox{ran}
V^{l+1}
\\ & = (V^l \clh)^\perp \oplus V^{l+1} \clh = \mbox{ker} V^{*l} \oplus V^{l+1} \clh,
\end{split}\]for all $l \geq 0$. Consequently, we have \[\clh_u =
\mathop{\bigcap}_{m=0}^\infty (\mbox{ker} V^{*m} \oplus V^{m+1}
\clh\big) = \mathop{\bigcap}_{m=0}^\infty V^m \clh.\] Uniqueness of
the decomposition readily follows from the uniqueness of the
wandering subspace $\clw$ for $V$. This completes the proof. \qed

We now introduce some notation which will remain fixed for the rest
of the paper. Given an integer $1 \leq m \leq n$, we denote the set
$\{1, \ldots, m\}$ by $I_m$. In particular, $I_n = \{1, \ldots,
n\}$.

Define $\clw_i := \mbox{ran}(I - V_i V_i^*)$ for each $1 \leq i \leq
n$ and
\[\clw_A := \mbox{ran}(\mathop{\prod}_{i \in A}
(I - V_i V_i^*)),\]where $A$ is a non-empty subset of $I_m$ and $1
\leq m \leq n$.

By doubly commutativity of $V$ it follows that
\[(I - V_i V_i^*) (I - V_j V_j^*) = (I - V_j V_j^*) (I - V_i
V_i^*),\]for all $i \neq j$. In particular, \[\{(I -V_i
V_i^*)\}_{i=1}^n\]is a family of commuting orthogonal projections on
$\clh$. Therefore, for all non-empty subsets $A$ of $I_m$ ($1 \leq m
\leq n$) we have
\begin{equation}\label{w-alpha}\clw_A = \mbox{ran}(\mathop{\prod}_{i \in A}
(I - V_i V_i^*)) = \mathop{\bigcap}_{i \in A} \mbox{ran} (I - V_i
V_i^*)) = \mathop{\bigcap}_{i \in A} \clw_i.\end{equation}

The following simple result plays a basic role in describing the
class of $n$-tuples of doubly commuting isometries.

\begin{Proposition}\label{A-j reducing}
Let $V = (V_1, \ldots, V_n)$ be an $n$-tuple ($n \geq 2$) of doubly
commuting isometries on $\clh$ and $A$ be a non-empty subset of
$I_m$ for $1 \leq m \leq n$. Then $\clw_A$ is a $V_j$-reducing
subspace of $\clh$ for all $j \in I_n \setminus A$.
\end{Proposition}

\NI\textsf{Proof.} By doubly commutativity of $V$ we have
\[V_j (I - V_i V_i^*) = (I - V_i V_i^*) V_j,\]for all $i \neq j$, and thus \[V_j (\mathop{\prod}_{i \in A} (I
- V_i V_i^*)) = (\mathop{\prod}_{i \in A} (I - V_i V_i^*)) V_j,\]for
all $j \in I_n \setminus A$, that is,
\[V_j P_{\clw_A} = P_{\clw_A} V_j,\]where $P_{\clw_A}$ is the orthogonal projection of $\clw_A$ onto $\clh$. This completes the proof. \qed

To complete this section we will use the preceding proposition to
obtain the generalized wandering subspaces for $n$-tuple of doubly
commuting isometries.

\begin{Corollary}\label{alpha-j}
Let $V = (V_1, \ldots, V_n)$ be an $n$-tuple ($n \geq 2$) of doubly
commuting isometries on $\clh$ and $m \leq n$. Then for each
non-empty subset $A$ of $I_m$ and $j \in I_n \setminus A$,
\[\clw_A \ominus V_j \clw_A = \mbox{ran}
(\mathop{\prod}_{i \in A} (I - V_i V_i^*) (I - V_j V_j^*)) =
\Big(\mathop{\bigcap}_{i \in A} \clw_i\Big)\; \cap \clw_j.\]
\end{Corollary}

\NI\textsf{Proof.} Doubly commutativity of $V$ implies that
\[\mathop{\prod}_{i \in A} (I - V_i V_i^*) (I - V_j V_j^*) =
\mathop{\prod}_{i \in A} (I - V_i V_i^*) - V_j (\mathop{\prod}_{i
\in A} (I - V_i V_i^*)) V_j^*.\]By Proposition \ref{A-j reducing} we
have $V_j \clw_A \subseteq \clw_A$ for all $j \not \in A$. Moreover
\[V_j \clw_A= \mbox{ran} [ V_j \mathop{\prod}_{i \in A} (I - V_i
V_i^*) V_j^*],\]and hence
\[\clw_A \ominus V_j \clw_A = \mbox{ran}
(\mathop{\prod}_{i \in A} (I - V_i V_i^*) - V_j (\mathop{\prod}_{i
\in A} (I - V_i V_i^*)) V_j^*) = \mbox{ran} (\mathop{\prod}_{i \in
A} (I - V_i V_i^*) (I - V_j V_j^*)),\] for all $j \not \in A$. The
second equality follows from (\ref{w-alpha}). This completes the
proof. \qed

\newsection{The Main Theorem}

In this section we will prove the main result of this paper.

Before proceeding, we shall adopt the following set of notations.
Let $(T_1, \ldots, T_n)$ be an $n$-tuple of commuting operators on a
Hilbert space $\clh$ and $1 \leq m \leq n$. Let $A = \{i_1, \ldots,
i_l\} \subseteq I_m$ and $1 \leq l \leq m$. We denote by $T_A$ the
$|A|$-tuple of commuting operators $(T_{i_1}, \ldots, T_{i_l})$ and
$\mathbb{N}^A := \{\bm{k} = (k_{i_1}, \ldots, k_{i_l}): k_{i_j} \in
\mathbb{N}, 1 \leq j \leq l\}$. We also denote $T_{i_1}^{k_{i_1}}
\cdots T_{i_l}^{k_{i_l}}$ by $T_A^{\bm{k}}$ for all $\bm{k} \in
\mathbb{N}^{A}$.

\begin{Theorem}\label{dc-w} Let $V = (V_1, \ldots, V_n)$ be an
$n$-tuple ($n \geq 2$) of doubly commuting isometries on $\clh$ and
$m \in \{2, \ldots, n\}$. Then there exists $2^m$ joint $(V_1,
\ldots, V_m)$-reducing subspaces $\{\clh_A : A \subseteq I_m\}$
(counting the trivial subspace $\{0\}$) such
that\begin{equation}\label{W1}\clh = \mathop{\bigoplus}_{A \subseteq
I_m} \clh_A,\end{equation}where for each $A \subseteq I_m$,
\begin{equation}\label{W2} \clh_A = \mathop{\bigoplus}_{\bm{k} \in \mathbb{N}^A}
V_A^{\bm{k}} \Big(\mathop{\bigcap}_{\bm{j} \in \mathbb{N}^{I_m
\setminus A}} V_{I_m \setminus A}^{\bm{j}}
\clw_A\Big).\end{equation}In particular, there exists $2^n$
orthogonal joint $V$-reducing subspaces
 $\{\clh_{A}: A \subseteq
I_n\}$ such that
\[\clh = \mathop{\sum}_{A \subseteq  I_n} \oplus \clh_A,\]and for each $A \subseteq I_n$ and $\clh_A \neq \{0\}$, $V_i|_{\clh_A}$ is a
shift if $i \in A$ and unitary if $i \in I_n \setminus A$ for all $i
= 1, \ldots, n$. Moreover, the above decomposition is unique, in the
sense that \[\clh_A = \mathop{\bigoplus}_{\bm{k} \in \mathbb{N}^A}
V_A^{\bm{k}} \Big(\mathop{\bigcap}_{\bm{j} \in \mathbb{N}^{I_n
\setminus A}} V_{I_n \setminus A}^{\bm{j}} \clw_A\Big),\]for all $A
\subseteq  I_n$.
\end{Theorem}
\NI\textsf{Proof.} We shall prove the statement using mathematical
induction.

\NI For $m = 2$: By applying the Wold decomposition theorem, Theorem
\ref{Wold}, to the isometry $V_1$ on $\clh$ we have
\[\clh = \mathop{\bigoplus}_{k_1 \in \mathbb{N}} V_1^{k_1}
\clw_1 \mathop{\bigoplus} \Big(\mathop{\bigcap}_{k_1 \in \mathbb{N}}
V_1^{k_1} \clh.\Big)\]As $\clw_1$ is a $V_2$-reducing subspace, it
follows from the Wold decomposition theorem for the isometry
$V_2|_{\clw_1} \in \cll(\clw_1)$ that
\[\clw_1 = \mathop{\bigoplus}_{k_2 \in \mathbb{N}} V_2^{k_2} (\clw_1 \ominus V_2 \clw_1) \bigoplus \Big(\mathop{\bigcap}_{k_2 \in
\mathbb{N}} V_2^{k_2} \clw_1\Big) = \mathop{\bigoplus}_{k_2 \in
\mathbb{N}} V_2^{k_2} (\clw_1 \cap \clw_2) \Big(\bigoplus
\mathop{\bigcap}_{k_2 \in \mathbb{N}} V_2^{k_2} \clw_1\Big),\]where
the second equality follows from Corollary \ref{alpha-j}.
Consequently,
\[\begin{split}\clh & = \mathop{\bigoplus}_{k_1 \in \mathbb{N}}
V_1^{k_1} \clw_1 \bigoplus \Big(\mathop{\bigcap}_{k_1 \in
\mathbb{N}} V_1^{k_1} \clh\Big)\\ & = \mathop{\bigoplus}_{k_1 \in
\mathbb{N}} V_1^{k_1} \Big(\mathop{\bigoplus}_{k_2 \in \mathbb{N}}
V_2^{k_2} \Big(\clw_1 \cap \clw_2\Big) \bigoplus
\mathop{\bigcap}_{k_2 \in \mathbb{N}} V_2^{k_2} \clw_1\Big)
\bigoplus \mathop{\bigcap}_{k_1 \in \mathbb{N}} V_1^{k_1} \clh \\ &
= \mathop{\bigoplus}_{k_1, k_2 \in \mathbb{N}} V_1^{k_1} V_2^{k_2}
\Big(\clw_1 \cap \clw_2\Big) \mathop{\bigoplus}_{k_1 \in \mathbb{N}}
V_1^{k_1} \Big(\mathop{\bigcap}_{k_2 \in \mathbb{N}} V_2^{k_2}
\clw_1 \Big)\bigoplus \mathop{\bigcap}_{k_1 \in \mathbb{N}}
V_1^{k_1} \clh.
\end{split}\]
Furthermore, the Wold decomposition of the isometry $V_2$ on $\clh$
\[\clh = \mathop{\bigoplus}_{k_2 \in \mathbb{N}} V_2^{k_2} \clw_2
\bigoplus \mathop{\bigcap}_{k_2 \in \mathbb{N}} V_2^{k_2} \clh,\]
yields \[V_1^{k_1} \clh = \mathop{\bigoplus}_{k_2 \in \mathbb{N}}
V_2^{k_2} V_1^{k_1} \clw_2 \bigoplus \mathop{\bigcap}_{k_2 \in
\mathbb{N}} V_1^{k_1}V_2^{k_2} \clh,\]for all $k_1 \in \mathbb{N}$.
From this we infer that
\[\mathop{\bigoplus}_{k_1 \in \mathbb{N}} V_1^{k_1} \clh =
\mathop{\bigoplus}_{k_2 \in \mathbb{N}} V_2^{k_2} \Big(
\mathop{\bigcap}_{k_1 \in \mathbb{N}} V_1^{k_1} \clw_2 \Big)
\bigoplus \mathop{\bigcap}_{k_1, k_2 \in \mathbb{N}} V_1^{k_1}
V_2^{k_2} \clh.\]Therefore, \begin{equation}\label{n=2case} \clh =
\mathop{\bigoplus}_{\bm{k} \in \mathbb{N}^2} V^{\bm{k}} \Big(\clw_1
\cap \clw_2\Big) \mathop{\bigoplus}_{k_1 \in \mathbb{N}} V_1^{k_1}
\Big(\mathop{\bigcap}_{k_2 \in \mathbb{N}} V_2^{k_2} \clw_1 \Big)
\mathop{\bigoplus}_{k_2 \in \mathbb{N}} V_2^{k_2} \Big(
\mathop{\bigcap}_{k_1 \in \mathbb{N}} V_1^{k_1} \clw_2 \Big)
\bigoplus \mathop{\bigcap}_{\bm{k} \in \mathbb{N}^2} V^{\bm{k}}
\clh,\end{equation}that is, \[\clh = \mathop{\bigoplus}_{A \subseteq
I_2} \clh_A,\]where $\clh_A$ is as in (\ref{W2}) for each $A
\subseteq I_2$.

\NI For $m +1 \leq n$: Now let for $m < n$, we have $\clh =
\mathop{\bigoplus}_{A \subseteq I_m} \clh_{A}$, where for each
non-empty subset $A$ of $I_m$
\[\clh_A = \mathop{\bigoplus}_{\bm{k} \in \mathbb{N}^A}
V_A^{\bm{k}} \Big(\mathop{\bigcap}_{\bm{j} \in \mathbb{N}^{I_m
\setminus A}} V_{I_m \setminus A}^{\bm{j}} \clw_A\Big),\]and for $A
= \phi \subseteq I_m$,
\[\clh_{A} = \mathop{\bigcap}_{\bm{k} \in \mathbb{N}^m}
V_{I_m}^{\bm{k}} \clh.\] We claim that \[\clh = \mathop{\sum}_{A
\subseteq I_{m+1}} \oplus \clh_A.\] Since $\clw_A$ is
$V_{m+1}$-reducing subspace for all none-empty $A \subseteq I_m$, we
have
\[\begin{split}\clw_A & = \mathop{\bigoplus}_{k_{m+1} \in
\mathbb{N}} V_{m+1}^{k_{m+1}} (\clw_A \ominus V_{m+1} \clw_A)
\bigoplus \Big(\mathop{\bigcap}_{k_{m+1} \in \mathbb{N}} V_{m+1}^{k_{m+1}} \clw_A\Big)\\
& = \mathop{\bigoplus}_{k_{m+1} \in \mathbb{N}} V_{m+1}^{k_{m+1}}
(\mathop{\bigcap}_{i \in A} \clw_i\; \bigcap \clw_{m+1}) \bigoplus
\Big(\mathop{\bigcap}_{k_{m+1} \in \mathbb{N}} V_{m+1}^{k_{m+1}}
\clw_A\Big),\end{split}\] and hence it follows that
\[\begin{split} \clh_A & = \mathop{\bigoplus}_{\bm{k} \in \mathbb{N}^A}
V_A^{\bm{k}} \Big(\mathop{\bigcap}_{\bm{j} \in \mathbb{N}^{I_m
\setminus A}} V_{I_m \setminus A}^{\bm{j}} \clw_A\Big)\\
& = \mathop{\bigoplus}_{\bm{k} \in \mathbb{N}^A} V_A^{\bm{k}}
\Big(\mathop{\bigcap}_{\bm{j} \in \mathbb{N}^{I_m \setminus A}}
V_{I_m \setminus A}^{\bm{j}} \Big(\mathop{\bigoplus}_{k_{m+1} \in
\mathbb{N}} V_{m+1}^{k_{m+1}} (\mathop{\bigcap}_{i \in A} \clw_i\;
\cap \clw_{m+1}) \bigoplus \Big(\mathop{\bigcap}_{k_{m+1} \in
\mathbb{N}} V_{m+1}^{k_{m+1}} \clw_A\Big) \Big) \Big) \\ & =
\Big[\mathop{\mathop{\bigoplus}_{\bm{k} \in \mathbb{N}^A}}_{ k_{m+1}
\in \mathbb{N}}
 V_A^{\bm{k}}  V_{m+1}^{k_{m+1}}
\Big(\mathop{\bigcap}_{\bm{j} \in \mathbb{N}^{I_m \setminus A}}
V_{I_m \setminus A}^{\bm{j}} (\mathop{\bigcap}_{i \in A \cup\{m+1\}}
\clw_i)\Big)\Big] \bigoplus \Big[ \mathop{\bigoplus}_{\bm{k} \in
\mathbb{N}^A} V_A^{\bm{k}} \mathop{\bigcap}_{k_{m+1} \in \mathbb{N}}
V_{m+1}^{k_{m+1}} \clw_A\Big]
\end{split}\]

\NI Applying again the Wold decomposition to the isometry $V_{m+1}$
on $\clh$, we have
\[\clh = \mathop{\bigoplus}_{k_{m+1} \in \mathbb{N}}
V_{m+1}^{k_{m+1}} \clw_{m+1} \bigoplus \mathop{\bigcap}_{k_{m+1} \in
\mathbb{N}} V_{m+1}^{k_{m+1}} \clh,\] and hence for $A = \phi
\subseteq I_m$,
\[\begin{split}\clh_A & = \mathop{\bigcap}_{\bm{k}
\in \mathbb{N}^m} V_1^{k_1} \cdots V_m^{k_m} \clh \\ & =
\mathop{\bigcap}_{\bm{k} \in \mathbb{N}^m} V_1^{k_1} \cdots
V_m^{k_m} ( \mathop{\bigoplus}_{k_{m+1} \in \mathbb{N}}
V_{m+1}^{k_{m+1}} \clw_{m+1} \bigoplus
\mathop{\bigcap}_{k_{m+1} \in \mathbb{N}} V_{m+1}^{k_{m+1}} \clh)\\
& = \Big[\mathop{\bigoplus}_{k_{m+1} \in \mathbb{N}}
V_{m+1}^{k_{m+1}} (\mathop{\bigcap}_{\bm{k} \in \mathbb{N}^m}
V_1^{k_1} \cdots V_m^{k_m} \clw_{m+1})\Big] \bigoplus
\Big[\mathop{\bigcap}_{\bm{k} \in \mathbb{N}^{A \cup \{m+1\}}}
V_1^{k_1} \cdots V_m^{k_m} V_{m+1}^{k_{m+1}} \clh\Big].
\end{split} \]Consequently, \[\clh = \mathop{\sum \oplus}_{A
\subseteq I_{m+1}} \clh_A.\]It follows immediately from the above
orthogonal decomposition of $\clh$ that $V_i|_{\clh_A}$ is a shift
for all $i \in A$ and unitary for all $i \in I_n \setminus A$.

\NI The uniqueness part follows from the uniqueness of the classical
Wold decomposition of isometries and the canonical construction of
the present orthogonal decomposition. This completes the proof. \qed

Note that if $n = 2$, then (\ref{n=2case}) yields a new proof of
Slocinski's Wold-type decomposition of a pair of doubly commuting
isometries.

The following corollary is an $n$-variables analogue of the
wandering subspace representations of pure isometries (that is,
shift operators) on Hilbert spaces.

\begin{Corollary}\label{pure}
Let $V = (V_1, \ldots, V_n)$ be an $n$-tuple ($n \geq 2$) of doubly
commuting shift operators on $\clh$. Then \[\clw =
\mathop{\bigcap}_{i=1}^n \mbox{ran} (I - V_i V_i^*),\]is a wandering
subspace for $V$ and
\[\clh = \mathop{\sum \oplus}_{\bm{k} \in \mathbb{N}^n} V^{\bm{k}}
\clw.\]
\end{Corollary}
\NI\textsf{Proof.} Let us note first that the given condition is
equivalent to \[\mathop{\bigcap}_{k \in \mathbb{N}} V_i^k \clh =
\{0\},\]for all $1 \leq i \leq n$. Then the result readily follows
from the proof of Theorem \ref{dc-w}. \qed

Recall that a pair of $n$-tuples $V = (V_1, \ldots, V_n)$ and $W =
(W_1, \ldots, W_n)$ on $\clh$ and $\clk$, respectively, are said to
be unitarily equivalent if there exists a unitary map $U : \clh
\raro \clk$ such that $U V_i = W_i U$ for all $1 \leq i \leq n$.

The following corollary is a generalization of Theorem 1 in
\cite{S}.

\begin{Theorem}
Let $V = (V_1, \ldots, V_n)$ be an $n$-tuple ($n \geq 2$) of
commuting isometries on $\clh$. Then the following conditions are
equivalent:

(i) There exists a wandering subspace $\clw$ for $V$ such that
\[\clh = \mathop{\bigoplus}_{\bm{k} \in \mathbb{N}^n}  V^{\bm{k}}
\clw.\]

(ii) $V_m$ is a shift for all $m = 1, \ldots, n$ and $V$ is doubly
commuting tuple.

(iii) There exists $m \in \{1, \ldots, n\}$ such that $V_m$ is a
shift and the wandering subspace for $V_m$ is given by \[\clw_m =
\mathop{\bigoplus}_{{\bm{k} \in \mathbb{N}^n},\, k_m = 0}
\;V^{\bm{k}} (\mathop{\bigcap}_{i = 0}^n \clw_i).\]

(iv) $\clw: = \bigcap_{i=1}^n \clw_i$ is a wandering subspace for
$V$ and that $\clh = \mathop{\bigoplus}_{\bm{k} \in \mathbb{N}^n}
V^{\bm{k}} \clw$.

(v) $V$ is unitarily equivalent to $M_z = (M_{z_1}, \ldots,
M_{z_n})$ on $H^2_{\cle}(\mathbb{D}^n)$ for some Hilbert space
$\cle$ with $\mbox{dim}\, \cle = \mbox{dim}\, \clw$.
\end{Theorem}
\NI\textsf{Proof.} (i) implies (ii) : That $V_m$ is a shift, for all
$1 \leq m \leq n$, follows from the fact that \[\clh =
\mathop{\bigoplus}_{k \in \mathbb{N}} V_m^{k}
\Big(\mathop{\bigoplus}_{\bm{k} \in \mathbb{N}^n,\, k_m = 0}
V^{\bm{k}} \clw\Big).\]Now let $h \in \clh$ and that \[h =
\mathop{\sum}_{i=0}^\infty V_m^i f_i. \quad \quad \quad (f_i \in
\mathop{\bigoplus}_{\bm{k} \in \mathbb{N}^n,\, k_m = 0} V^{\bm{k}}
\clw)\]It follows that for all $l \neq m$, \[V_l V_m^* h = V_l
(\mathop{\sum}_{i=1}^\infty V_m^{i-1} f_i) =
\mathop{\sum}_{i=1}^\infty V_m^{i-1} (V_l f_i) = V_m^*
(\mathop{\sum}_{i=0}^\infty V_m^i (V_l f_i)) = V_m^* V_l
(\mathop{\sum}_{i=0}^\infty V_m^{i} f_i),\]that is, $V$ is doubly
commuting.

(ii) implies (iii): By Corollary \ref{dc-w} we have \[\clh =
\mathop{\bigoplus}_{\bm{k} \in \mathbb{N}^n} V^{\bm{k}}
(\mathop{\bigcap}_{i=1}^n \clw_i) = \mathop{\bigoplus}_{k \in
\mathbb{N}} V_m^{k} \Big(\mathop{\bigoplus}_{\bm{k} \in
\mathbb{N}^n,\, k_m = 0} V^{\bm{k}} (\mathop{\bigcap}_{i=1}^n
\clw_i) \Big),\] and hence (iii) follows.

(iii) implies (iv): Since $V_m$ is a shift with the wandering
subspace \[\clw_m = \mathop{\bigoplus}_{{\bm{k} \in \mathbb{N}^n},\,
k_m = 0} V^{\bm{k}} (\mathop{\bigcap}_{i = 0}^n \clw_i),\]we infer
that
\[\clh = \mathop{\bigoplus}_{k \in \mathbb{N}} V_m^k \clw_m =
\mathop{\bigoplus}_{k \in \mathbb{N}} V_m^k
\Big(\mathop{\bigoplus}_{{\bm{k} \in \mathbb{N}^n},\, k_m = 0}
V^{\bm{k}} (\mathop{\bigcap}_{i = 0}^n \clw_i)\Big) =
\mathop{\bigoplus}_{k \in \mathbb{N}} V_m^k
(\mathop{\bigcap}_{i=1}^n\clw_i),
\]and hence (iv) follows.

(iv) implies (v): Let $\cle = \mathop{\bigcap}_{i=1}^n\clw_i$.
Define the unitary operator \[U : \clh = \mathop{\bigoplus}_{\bm{k}
\in \mathbb{N}^n} V^{\bm{k}} (\mathop{\bigcap}_{i=1}^n\clw_i)
\longrightarrow H^2_{\cle}(\mathbb{D}^n) =
\mathop{\bigoplus}_{\bm{k} \in \mathbb{N}^n} z^{\bm{k}} \cle,\]by
$U(V^{\bm{k}} \eta) = z^{\bm{k}} \eta$ for all $\eta \in \cle$ and
$\bm{k} \in \mathbb{N}^n$. It is obvious that $U V_i = M_{z_i} U$
for all $i = 1, \ldots, n$.

That (v) implies (i) is trivial.

\NI This concludes the proof of the theorem. \qed

As we indicated earlier, the Wold decomposition theorem may be
regarded as a very powerful tool from which a number of classical
results may be deduced. In \cite{SSW} we will use techniques
developed in this paper to obtain a complete classification of
doubly commuting invariant subspaces of the Hardy space over the
polydisc.

It is worth mentioning that in \cite{P}, Popovici obtained a Wold
type decomposition for general isometries. See \cite{BDF}
(Proposition 2.8 and Corollary 2.9) and \cite{GG} for more results
along this line.

Finally, we point out that all results of this paper are also valid
for a family (that is, not necessarily finite tuple) of doubly
commuting isometries.

\vspace{0.2in}

\NI \textsf{Acknowledgement:} We would like to thank B. V. Rajarama
Bhat for valuable discussions and Joseph A. Ball for many useful
suggestions that improved the contents of the paper. Also we thank
Baruch Solel and Adam Skalski for pointing out the reference
\cite{SZ}.


\end{document}